# ON THE STABILITY OF GEOMETRIC EXTREMES

## S. Satheesh and N. Unnikrishnan Nair

*Cochin University of Science & Technology, Cochin.*

## Abstract

Possible reasons for the uniqueness of the geometric law in the context of stability of random extremes are explored here culminating in a conjecture.

## 1 Introduction

Let $F(x)$, $x \in \boldsymbol{R}$ be the distribution function (d.f) of a continuous random variable (r.v) $X$ and $Q(s)$ be the probability generating function (PGF) of a positive integer valued random variable (r.v) $N$ independent of $X$. Then $F(x)$ is N-max stable if

$$Q[F(x)] \;=\; F(a+bx), \tag{1}$$

and N-min stable if

$$Q[\,\overline{F}(a+bx)\,] = \overline{F}(x), \tag{2}$$

for all $x \in \boldsymbol{R}$ and some $a \in \boldsymbol{R}$ and $b>0$, where $\overline{F}(x) = 1\text{-}F(x)$. That is, the distributions of N-max and N-min of $F(x)$ should be of the same type as that of $F(x)$. For brevity we put $F_t(x) = F(a+bx)$ and $\overline{F}_t(x) \;=\; \overline{F}(a+bx)$.

The geometric law on {1,2,3, ....} with expectation $1/p$ will be denoted by geometric($p$). The PGF of the geometric($p$) law is $ps/[1-(1-p)s]$.

The semi-Pareto family of laws was characterized, among continuous distributions on [0,∞), by geometric($p$)-max stability in Pillai (1991) and by geometric($p$)-min stability in Pillai and Sandhya (1996). From these two characterizations it is clear that among distributions with non-negative support, geometric($p$)-max stability implies geometric($p$)-min stability and vice-versa, as both identify the same family. A natural curiosity thus is whether we can prove this without referring to the family of semi Pareto laws and also whether it is true in general for d.fs with support ***R.*** Proving that this indeed is true leads to the question whether it is unique of the geometric($p$) law. Unearthing certain properties of the geometric($p$) PGF implicit in this proof, we obtain necessary and sufficient conditions on the PGF of $N$ for the stability of N-extremes (both N-max and N-min). In this attempt to characterize the geometric($p$) law we arrive at a conjecture, and this is done in section 2.

Sreehari (1995) has characterized the geometric($p$) law by N-max-stability of the semi-Pareto family while Satheesh and Nair (2002) by N-min-stability. The present study is also significant in the following two (similar) contexts. Arnold, *et al*. (1986) made the remarkable observation that the distribution of geometric($p$)-min of geometric($p$)-maxs of $F(x)$ has the same algebraic structure as that of geometric($p$)-max of geometric($p$)-mins. Marshall and Olkin (1997) introduced a parameterization scheme for a survival function $\overline{F}(x)$, $x \in \boldsymbol{R}$ that is similar in structure to the geometric($p$)-minimums by defining another survival function

$$\overline{G}(x,a) \quad = \frac{a\overline{F}(x)}{1-(1-a)\overline{F}(x)}, x \in \boldsymbol{R}, \ a > 0 \tag{3}$$

and showed that this family is geometric($p$)-extreme (that is, both geometric($p$)-max and geometric($p$)-min) stable. They attributed this property partially to the fact that geometric($p$) laws are closed under their own compounding. They also concluded that the N-min stability could not be expected if the geometric($p$) distribution is replaced by another distribution. See Remark.2.7 in this context. The observations made by Arnold, *et al*. (1986) and Marshall and Olkin (1997) are closely related.

In the concluding remarks we discuss these two situations and give the reason behind this property and supplement the arguments of Marshall and Olkin. Also, our conjecture comes closer in justifying the uniqueness of the geometric($p$) law in the context.

## 2 Uniqueness of the Geometric($p$) law

**Theorem.2.1** A d.f $F(x)$, $x \in \boldsymbol{R}$ is geometric($p$)-max stable if and only if it is geometric($p$)-min stable. (Notice that the geometric($p$) laws have the same parameter.)

**Proof.** We have:

$$\frac{pF(x)}{\left[1-(1-p)F(x)\right]} = F_t(x) \tag{1}$$

$$\Leftrightarrow \quad 1 - \frac{pF(x)}{\left[1-(1-p)F(x)\right]} = \overline{F}_t(x) \tag{2}$$

$$\Leftrightarrow \quad \frac{\overline{F}(x)}{p+(1-p)\overline{F}(x)} = \overline{F}_t(x) \tag{3}$$

$$\Leftrightarrow \quad \overline{F}(x) = \frac{p\overline{F}_t(x)}{1-(1-p)\overline{F}_t(x)} \tag{4}$$

Here (2.1) represents geometric($p$)-max stability and (2.4) geometric($p$)-min stability of $F(x)$ and the proof is complete.

The curiosity now is whether N-max and N-min stability of $F(x)$ with respect to the same $N$ implies $N$ is geometric($p$). A closer look at the above proof in this context reveals certain properties of $Q_p(s)$, the PGF of the geometric($p$) law. The L.H.S of (2.2) describes the survival function of

geometric($p$)-max of $F(x)$. Writing this survival function in terms of $\overline{F}(x)$ and equating to the L.H.S of (2.3) shows that $Q_p(s)$ satisfies

$$1^- Q_p(1\text{-}s) = Q_p^{-1}(s). \tag{5}$$

From (2.3) a subsequent inversion resulted in the geometric($p$)-min stability of $F(x)$ in (2.4).

The L.H.S of (2.3) which is $Q_p^{-1}[\overline{F}(x)]$, specifies the distribution (survival function) of the geometric($p$)-max of $F(x)$ in terms of $\overline{F}(x)$. This can also be written as:

$$\frac{\lambda\overline{F}(x)}{1-(1-\lambda)\overline{F}(x)}, \lambda=1/p. \tag{6}$$

**Remark.2.1** Thus the survival functions of geometric($q$)-min and geometric($p$)-max of $F(x)$ have the same algebraic structure in terms of $\overline{F}(x)$.

Now taking geometric($q$)-min of these geometric($p$)-maxs (*ie*.(2.6)), the resulting distribution is specified by the survival function

$$\frac{q\lambda\overline{F}(x)}{1-(1-q\lambda)\overline{F}(x)}. \tag{7}$$

Notice that this survival function also has the same algebraic structure in terms of $\overline{F}(x)$ as that of the geometric($q$)-min. This has been possible because of the following two facts. (i) The L.H.S of (2.3) could be written as (2.6) and (ii) the PGFs of independent geometric($p$) laws are closed under their own compounding. That is,

$$Q_p^{-1}(s) = Q_\lambda(s), \lambda = 1/p \text{ and} \tag{8}$$

$$Q_p[Q_q(s)] = Q_{pq}(s) \quad \text{for all } |s|<1 \text{ and } 0<p,q<1. \tag{9}$$

**Remark.2.2** Here (2.8) and (2.5) imply that the survival functions of geometric($p$)-max and geometric($p$)-min of $F(x)$ have the same algebraic structure in terms of $\overline{F}(x)$.

**Remark.2.3** As a consequence of (2.9) and (2.5), the structure of (2.6) is retained even if we take geometric($q$)-max instead of geometric($q$)-min of (2.6). This is because we can retrace (2.6), (2.3), (2.2) and (2.1), take the geometric($q$)-max of (2.1) and come to the form of (2.6).

Motivated by these observations on the PGF of the geometric($p$) law we prove the following general results. Here $Q_u(s)$ denotes the PGF of $N$ with parameter $u>0$.

**Theorem.2.2** N-max stability of $F(x)$ implies N-min stability of $F(x)$ (and vice-versa) if and only if $Q_u(s)$ satisfies

$$1 - Q_u(1-s) = Q_u^{-1}(s), \, \forall \, 0<s<1. \tag{5b}$$

**Proof.** Suppose $F(x)$ is N-max stable and $Q_u(s)$ satisfies (2.5$b$). Then:

$$Q_u[F(x)] = F_t(x)$$

$$\Leftrightarrow \quad 1 - Q_u[1 - \overline{F}(x)] = \overline{F}_t(x)$$

$$\Leftrightarrow \quad Q_u^{-1}[\overline{F}(x)] = \overline{F}_t(x) \quad \text{(by (2.5$b$))}$$

$$\Leftrightarrow \quad \overline{F}(x) = Q_u[\overline{F}_t(x)] \text{ and hence } F(x) \text{ is N-min stable.}$$

Conversely, suppose N-max stability of $F(x)$ implies N-min stability of $F(x)$. Then we have, in terms of survival functions:

$$1 - Q_u[1 - \overline{F}(x)] = \overline{F}_t(x) \implies \overline{F}(x) = Q_u[\overline{F}_t(x)]$$

$$\Leftrightarrow Q_u^{-1}[\overline{F}(x)] = \overline{F}_t(x).$$

Hence,

$$1 - Q_u[1 - \overline{F}(x)] = Q_u^{-1}[\overline{F}(x)] \text{ or } 1 - Q_u(1-s) = Q_u^{-1}(s)$$

which is condition (2.5*b*) and the proof is complete.

**Theorem.2.3** Suppose the PGF $Q_u(s)$ of $N$ satisfies (2.5*b*). Then the survival functions of N-max and N-min of $F(x)$ have the same algebraic structure in terms of $\overline{F}(x)$ if and only if

$$Q_u^{-1}(s) = Q_\lambda(s), \text{ for some } \lambda > 0. \tag{8b}$$

**Proof.** The survival functions of N-max in terms of $\overline{F}(x)$ is $Q_u^{-1}[\overline{F}(x)]$ as $Q_u(s)$ satisfies (2.5*b*). Now the requirement implies that $Q_u^{-1}[\overline{F}(x)]$ and $Q_u[\overline{F}(x)]$ should have the same algebraic structure. Hence $Q_u^{-1}(s) = Q_\lambda(s)$ for some $\lambda > 0$. The converse is clear.                    €

**Theorem.2.4** Suppose the PGF $Q_u(s)$ of $N$ satisfies (2.5*b*) and (2.8*b*). If $Q_u(s)$ also satisfies

$$Q_u[Q_v(s)] = Q_{uv}(s) \quad \text{for all} \quad |s| < 1 \text{ and } u, v > 0, \tag{9b}$$

then the survival functions of N-min of N-maxs and N-max of N-mins of $F(x)$ have the same algebraic structure as that of N-min of $F(x)$.

**Proof.** The survival function of N-min of N-maxs of $F(x)$ is $Q_u\{Q_\lambda[\overline{F}(x)]\}$. If $Q_u(s)$ satisfies (2.9*b*) then

$$Q_u\{Q_\lambda[\overline{F}(x)]\} = Q_{u\lambda}[\overline{F}(x)].$$

This has the same structure as that of $Q_u[\overline{F}(x)]$. Similar is the case with N-max of N-mins. Hence the proof is complete.                    €

**Remark.2.4** Thus the reason for the phenomenon observed by Arnold, *et al*. (1986) is that the geometric($p$) PGF satisfies (2.5*b*), (2.8*b*) and (2.9*b*).

**Remark.2.5** From the proof of Theorem.2.4 it is also clear that when we consider the N-max of N-min (or reverse) of $F(x)$ it is not necessary that $N$'s have the same parameter at both places.

We now show that (2.9*b*) implies (2.8*b*).

**Lemma.2.1** If a one-to-one function $Q_u(s)$, $u>0$ satisfies (2.9*b*) then it satisfies (2.8*b*) with $\lambda = 1/u$ and $Q_1(s) = s$ for all $s$.

**Proof.** We have $Q_u[Q_v(s)] = Q_{uv}(s)$.

When $s = Q_{1/v}(s)$, $Q_u[Q_v[Q_{1/v}(s)]] = Q_u(s)$, which shows that

$Q_{1/v}(s) = Q_v^{-1}(s)$ and $Q_1(s) = s$ for all $s$.

**Remark.2.6** In Theorem.2.4, if $N$'s have the same parameter at both N-max and N-min then the survival function of N-max of N-min (or reverse) of $F(x)$ equals $\overline{F}(x)$ by Lemma.2.1.

Next we consider certain examples showing that there are PGFs that satisfy conditions (2.5*b*) or (2.9*b*) (hence (2.8*b*) also) other than the geometric($p$). Notice that equation (2.5*b*) can be equivalently written as:

$$Q_u[1 - Q_u(1 - s)] = s. \tag{10}$$

First we demonstrate that the conclusion of Theorem.2.1 is true for a non-geometric($p$) law for $N$. For instance, from an example in Shaked (1975) we have:

**Example.2.1** If $F(x)$ is N-max stable where $N$ has the PGF $1 - (1 - s^m)^{1/m}$, $m>1$ integer, then

$$1 - \{1 - [F(x)]^m\}^{1/m} = F_t(x)$$

$$\Leftrightarrow \qquad 1 - [F(x)]^m = [\overline{F}_t(x)]^m$$

$$\Leftrightarrow \qquad F(x) = \{1 - [\overline{F}_t(x)]^m\}^{1/m}$$

$$\Leftrightarrow \qquad \overline{F}(x) = 1 - \{1 - [\overline{F}_t(x)]^m\}^{1/m}.$$

Hence $F(x)$ is N-min stable as well. Clearly, the converse is also true.

Shaked (1975) arrived at (2.10) from the requirement that N-min of N-maxs of $F(x)$ must be stable (though he doesn't use these terms). He also solved the functional equation (2.10) (and thus (2.5*b*)) under the assumption that $Q(s)$ is single valued and meromorphic in the complex plane to characterize the geometric($p$) PGF. Now, can we restrict our search more realistically so as to characterize the geometric($p$) law?

The question is whether conditions (2.5*b*) and (2.9*b*) characterize the PGF of the geometric($p$) law. For example, the PGF of the Harris($a$,$k$) law

$$\frac{s}{\{a-(a-1)s^k\}^{1/k}}, \; k>0 \text{ integer and } a>1 \tag{11}$$

satisfies (2.9*b*). We may easily verify by direct computation using PGFs that the Harris($a$,$j$) laws are closed under their own compounding. We record this as:

**Lemma.2.2** Let $P_{\mathrm{u}}(s) = \dfrac{s}{[u-(u-1)s^j]^{1/j}}$ and $Q_{\mathrm{v}}(s) = \dfrac{s}{[v-(v-1)s^j]^{1/j}}$.

Then, $P_{\mathrm{u}}(Q_{\mathrm{v}}(s)) = \dfrac{s}{[uv-(uv-1)s^j]^{1/j}}$. Also $P_{\mathrm{u}}^{-1}(s) = P_\lambda(s)$, $\lambda = 1/u$ .

Notice that the parameter $j$ must be the same for both the PGFs and it is the same in the compound also. Thus it also satisfies (2.8*b*) with $\lambda = 1/u$, but it is not a solution of (2.5*b*). Also none of the examples in Shaked (1975) (including that in example.2.1) satisfies (2.9*b*) though they are solutions of (2.5*b*). Thus there are PGFs that satisfy either (2.5*b*) or (2.9*b*). Hence it appears that the only PGF that satisfies (2.5*b*) and (2.9*b*) is that of the geometric($p$) law. We state this as a conjecture.

**Conjecture.2.1** A PGF $Q_{\mathrm{u}}(s)$, $u>0$ satisfies (2.5*b*) and (2.9*b*) (and hence (2.8*b*)) if and only if it is the PGF of the geometric($u$) law with mean $1/u$.

One can generalize the Marshall and Olkin (1997) parameterization scheme on the following lines. This may be useful in lending more flexibility to the d.f. $F(x)$ in modeling. In Lemma.2.2 let $P_{\mathrm{u}}(s)$ be the PGF of a r.v. $N$ and $Q_{\mathrm{v}}(s)$ be that of $M$. Let $\{X_i\}$ be independent copies of a r.v.

$X$ with d.f $F(x)$, $x \in \boldsymbol{R}$. Let $N$, $M$ and $X$ be mutually independent. Put $U = Min(X_1, \ldots, X_N)$. Then,

$$P\{U > x\} = \frac{\overline{F}(x)}{[u - (u-1)[\overline{F}(x)]^j]^{1/j}}, x \in \boldsymbol{R}, \ j > 0 \text{ integer and } u > 1. \qquad (12)$$

**Proposition.2.1** The family of distributions of the form (2.12) is M-min stable.

**Proof.** Let $U_1$, $U_2$, $\ldots$ be independent copies of $U$ and $N_1$, $N_2$, $\ldots$ independent copies of $N$. Then (as in the proof of Proposition.5.1 in Marshall and Olkin (1997)) we have,

$$Min(U_1, \ldots, U_M) \ = Min(X_{11}, \ldots, X_{1N_1}, \ldots, X_{M1}, \ldots, X_{MN_M})$$

$$= Min(X_1, \ldots, X_{N_1 + \ldots + N_M})$$

by re-indexing $X_{ij}$. Now by virtue of Lemma.2.2, $N_1 + \ldots + N_M$ has a Harris($uv$,$j$) distribution. Hence the distribution of $Min(U_1, \ldots, U_M)$ is specified by a survival function of the form (2.12) with $uv$ instead of $u$. Hence the result is proved.

Again, setting $V = Max(X_1, \ldots, X_N)$ we have:

$$P\{V < x\} = \frac{F(x)}{\{u - (u-1)[F(x)]^j\}^{1/j}}, x \in \boldsymbol{R}, \ j > 0 \text{ integer and } u > 1. \qquad (13)$$

Now, proceeding as in the proof of Proposition.2.1 we have:

**Proposition.2.2** The family of distributions of the form (2.13) is M-max stable.

**Remark.2.7** Thus we do have N-min and N-max stability of $F(x)$ with respect to a non-geometric($p$) (Harris) r.v $N$. Though Harris laws are closed under their own compounding the families (2.12) and (2.13) are not Harris-extreme stable, because the PGF of the Harris law is not a solution of (2.5$b$).

**Concluding Remarks.** Arnold, *et al*. (1986) considered geometric($p$)-mins and geometric($p$)-maxs applied one after the other (in any order) and observed the similarity in the functional forms of their survival functions. The reason for this phenomenon is given in Remark.2.4. Marshall and Olkin (1997) has shown that their parameterization scheme, similar to the form in (2.6) with $\lambda$>0, is geometric($p$)-extreme stable. This is similar to that of Arnold, *et al*. (1986) though they are apparently different in form. They (Marshall and Olkin) had attributed this property partially to the fact that geometric($p$) laws are closed under their own compounding that is equivalent to (2.9). The reason given by Marshall and Olkin, which explains only part of the picture, is complemented by our observations (Remarks 2.2, 2.3 and 2.4) that the structure of (2.6) is retained because the geometric($p$) PGFs also satisfy (2.5$b$) in addition to (2.9). This resulted in Theorem.2.4 and subsequently led us to the Conjecture.2.1.

Finally, a complete proof of the conjecture will also show that among d.fs $F(x)$ with non-negative support, N-extreme stability will imply that $N$ is geometric($p$) and consequently $F(x)$ is semi Pareto. This will be a simultaneous characterization of the geometric($p$) and the semi-Pareto laws.

## Acknowledgement

The authors thank the referee for a careful reading of the paper and many helpful suggestions for clarity and improvement and to include remarks 2.5 and 2.6.

## References

**Arnold, B.C; Robertson, C.A; and Yeh, H.C.** (1986). Some properties of a Pareto type distribution, *Sankhya-A*, **48**, 404-408.

**Marshall, A.W; and Olkin, I.** (1997). On adding a parameter to a distribution with special reference to exponential and Weibull models, *Biometrika*, **84**, 641-652.

**Pillai, R.N.** (1991). Semi-Pareto processes, *J. Appl. Prob.*, **28**, 461-465.

**Pillai, R.N; and Sandhya, E.** (1996). Geometric sums and Pareto law in reliability theory, *I. A. P. Q. R. Trans.*, **21**, 137-142.

**Satheesh, S and Nair, N.U.** (2002). A note on maximum and minimum stability of certain distributions, *Calcutta Statist. Assoc. Bulletin*, **53**, 249-252.

**Shaked, M.** (1975). On the distributions of the minimum and of the maximum of a random number of i.i.d random variables, *In Statistical Distributions in Scientific Work*, Editors, Patil, G.P; Kotz, S; and Ord, J.K., D.Reidel Publishing Company, Dordrecht, Holland, 363-380.

**Sreehari, M.** (1995). Maximum stability and a generalization, *Statist. Probab. Lett.*, **23**, 339-342.

**S Satheesh**

Department of Statistics

Cochin University of Science and Technology

Cochin - 682 022, **India**.

ssatheesh@sancharnet.in

**N Unnikrishnan Nair**

Department of Statistics

Cochin University of Science and Technology

Cochin - 682 022, **India**.